\documentclass[12pt,reqno]{amsart}
\usepackage{amssymb}
\usepackage{amscd}
\usepackage{amsmath}

\newtheorem{theorem}{Theorem}[section]

\theoremstyle{definition}

\theoremstyle{remark}

\theoremstyle{coroll}
\newtheorem{coroll}[theorem]{Corollary}

\theoremstyle{conject}
\newtheorem{conject}[theorem]{Conjecture}

\numberwithin{equation}{section}

\newcommand{\nc}{\newcommand}
\nc{\cal}{\mathcal} 
\nc{\la}{\langle} \nc{\ra}{\rangle}
 \nc{\CA}{\cal A}
 \nc{\CBB}{\cal B}
 \nc{\CC}{\cal C}
\nc{\CDD}{\cal D} \nc{\CE}{\cal E} \nc{\CF}{\cal F} \nc{\CG}{\cal
G} \nc{\CH}{\cal H} \nc{\CI}{\cal I} \nc{\CJ}{\cal J}
\nc{\CK}{\cal K} \nc{\CL}{\cal L} \nc{\CM}{\cal M} \nc{\CN}{\cal
N} \nc{\CO}{\cal O} \nc{\CP}{\cal P} \nc{\CQ}{\cal Q}
\nc{\CR}{\cal R} \nc{\CS}{\cal S} \nc{\CT}{\cal T} \nc{\CU}{\cal
U} \nc{\CV}{\cal V} \nc{\CW}{\cal W} \nc{\CZ}{\cal Z}

\nc{\kap}{\kappa}

\nc{\fa}{\mathfrak a} \nc{\fg}{\mathfrak g} \nc{\fk}{\mathfrak k}
\nc{\fh}{\mathfrak h} \nc{\fm}{\mathfrak m} \nc{\fn}{\mathfrak n}
\nc{\fA}{\mathfrak A} \nc{\fC}{\mathfrak C} \nc{\fI}{\mathfrak I}
\nc{\fL}{\mathfrak L} \nc{\fS}{\mathfrak S}

\nc{\nen}{\newenvironment} \nc{\ol}{\overline}
\nc{\ul}{\underline} \nc{\lra}{\longrightarrow}
\nc{\lla}{\longleftarrow} \nc{\Lra}{\Longrightarrow}
\nc{\Lla}{\Longleftarrow} \nc{\Llra}{\Longleftrightarrow}
\nc{\hra}{\hookrightarrow} \nc{\iso}{\overset{\sim}{\lra}}
\nc{\Hom}{\mathrm{Hom}} \nc{\Mor}{\mathrm{Mor}}

\nc{\notebox}[1]{\noindent\fbox{\parbox{12.5cm}{\sf #1}}\\[8pt]}
\nc{\Thm}[1]{Theorem~\ref{#1}} \nc{\Prop}[1]{Proposition~\ref{#1}}
\nc{\Lem}[1]{Lemma~\ref{#1}} \nc{\Cor}[1]{Corollary~\ref{#1}}
\nc{\Conj}[1]{Conjecture~\ref{#1}} \nc{\Claim}[1]{Claim~\ref{#1}}
\nc{\Defn}[1]{ Definition~\ref{#1}} \nc{\Exa}[1]{Example~\ref{#1}}
\nc{\Rem}[1]{Remark~\ref{#1}} \nc{\Note}[1]{Note~\ref{#1}}

\marginparsep 0.1cm \marginparwidth 2.5cm \nc{\marg}{\marginpar}

\nen{thm}[1]{\label{#1}{\bf Theorem.\ } \em}{}
\nen{thma}[1]{\label{#1}{\bf Theorem A.\ } \em}{}
\nen{thmb}[1]{\label{#1}{\bf Theorem B.\ } \em}{}
\nen{prop}[1]{\label{#1}{\bf Proposition.\ } \em}{}
\nen{lem}[1]{\label{#1}{\bf Lemma.\ } \em}{}
\nen{klem}[1]{\label{#1}{\bf Key Lemma.\ } \em}{}
\nen{cor}[1]{\label{#1}{\bf Corollary.\ } \em}{}
\nen{cora}[1]{\label{#1}{\bf Corollary A.\ } \em}{}
\nen{corb}[1]{\label{#1}{\bf Corollary B.\ } \em}{}
\nen{conj}[1]{\label{#1}{\bf Conjecture.\ } \em}{}
\nen{claim}[1]{\label{#1}{\bf Claim.\ } \em}{}

\nen{defn}[1]{\label{#1}{\bf Definition.\ } }{}
\nen{exa}[1]{\label{#1}{\bf Example.\ } }{}

\nen{rem}[1]{\label{#1}{\em Remark.\ } }{}
\nen{exer}[1]{\label{#1}{\em Exercise.\ } }{}
\nen{pf}{\begin{proof}}{\end{proof}}

 \nc{\br}{\mathbb R}
 \nc{\bz}{\mathbb Z}
 \nc{\bc}{\mathbb C}
 \nc{\bn}{\mathbb N}
 \nc{\geg}{\mathfrak g}
 \nc{\G}{\Gamma}
 \nc{\sm}{\setminus}
 \nc{\sub}{\subset}
 \nc{\lm}{\lambda}
 \nc{\al}{\alpha}
 \nc{\bt}{\beta}
 \nc{\om}{\omega}
 \nc{\dl}{\delta}
 \nc{\g}{\gamma}
 \nc{\Dl}{\Delta}
 \nc{\Om}{\Omega}
 \nc{\s}{\sigma}
 \nc{\ro}{\rho}
 \nc{\te}{\theta}
 \nc{\SLR}{SL_2(\br)}
 \nc{\GLR}{GL_2(\br)}
 \nc{\PGLR}{PGL_2(\br)}
 \nc{\PSLR}{PSL_2(\br)}
 \nc{\SLC}{SL(2,\bc)}
 \nc{\uH}{\mathbb H}
 \nc{\fD}{\mathcal{D}}
 \nc{\fE}{\mathcal{E}}
 \nc{\fO}{\mathcal{O}}
 \nc{\haf}{\frac{1}{2}}
 \nc{\qtr}{\frac{1}{4}}
 \nc{\shaf}{{\scriptstyle\frac{1}{2}}}

 \nc{\8}{\infty}
 \nc{\7}{{-\infty}}
 \nc{\inv}{^{-1}}
 \nc{\eps}{\varepsilon}

\begin{document}

\title[periods and representation theory]
{Periods, Subconvexity of $L$-functions and Representation Theory}

\author{Joseph Bernstein }
\address{Tel Aviv University, Ramat Aviv, Israel}
\email{bernstei@post.tau.ac.il}


\author{Andre Reznikov}
\address{Bar Ilan University, Ramat-Gan, Israel}
\email{reznikov@math.biu.ac.il}


\date{\today}

\dedicatory{Dedicated to Raoul Bott.}


\begin{abstract} We describe a new method to estimate the
trilinear period on automorphic representations of $\PGLR$. Such a
period gives rise to a special value of the triple $L$-function.
We prove a bound for the triple period which amounts to a
subconvexity bound for the corresponding special value. Our method
is based on the study of the analytic structure of the
corresponding unique trilinear functional on unitary
representations of $\PGLR$.
\end{abstract}
\maketitle
\section{Introduction}
\label{intro}

\subsection{Maass forms} \ Let $\uH$ denote the upper half plane equipped with the
standard Riemannian metric of constant curvature $-1$. We denote by $dv$  the
associated volume element  and by $\Delta$ the corresponding Laplace-Beltrami
operator on $\uH$.

    Fix a discrete group $\G$ of motions of $\uH$ and consider the Riemann surface
 $Y = \G \backslash \uH$. For simplicity we assume that $Y$ is compact (the case of
 $Y$ of finite volume is discussed at the end of the introduction). According to
  the uniformization theorem any compact Riemann surface $Y$ with the metric of
  constant curvature $-1$ is a special case of this construction

    Consider the spectral decomposition of the operator $\Delta$ in the space
    $L^2(Y,dv)$ of functions on $Y$. It is known that the operator $\Delta$
    is non-negative and has purely discrete spectrum;
    we will denote by $0=\mu_0< \mu_1 \leq \mu_2 \leq ...$ the
eigenvalues of $\Delta$. For these eigenvalues we always use a
natural from representation-theoretic point of view
parametrization $\mu_i=\frac{1-\lm_i^2}{4}$, where $\lm_i\in \bc$.
We denote by $\phi_i=\phi_{\lm_i}$ the corresponding
eigenfunctions (normalized to have $L^2$-norm one).

In the theory of automorphic forms, the functions $\phi_{\lm_i}$
are called automorphic functions or {\it Maass forms} (after H.
Maass, \cite{M}). The study of Maass forms plays an important role
in analytic number theory, analysis and geometry. We are
interested in their analytic properties and will present a new
method of bounding some important quantities arising from the
$\phi_{i}$.

A particular problem we are going to address in this paper belongs
to an active area of research in the theory of automorphic
functions studying an interplay between periods, special values of
automorphic $L$-functions and representation theory. One of the
central features of this interplay is the uniqueness of invariant
functionals associated with corresponding periods. The discovery
of this interplay goes back to classical works of E. Hecke and H.
Maass.

It is well-known that the uniqueness principle plays a central
role in the modern theory of automorphic functions (see
\cite{PS}). The impact uniqueness has on the analytic behavior of
periods and $L$-functions is yet another manifestation of this
principle.

\subsection{Triple products}\label{triple} \  For any three Maass
forms $\phi_i,\ \phi_j, \ \phi_k$, we define the following {\it triple product} or
{\it triple period}:\label{cijk}
\begin{eqnarray}
c_{ijk}=\int_Y\phi_i\phi_j\phi_kdv \ .
\end{eqnarray}

We would like to estimate the coefficient $c_{ijk}$ as a function
of parameters $\lm_i,\ \lm_j,\ \lm_k$. In particular, we would
like to find bounds for these coefficients as one or more of the
indices $i,\ j,\ k$ tend to infinity.

The bounds on the coefficient $c_{ijk}$ are related to bounds on
automorphic $L$-functions as can be seen from the following
beautiful formula of T. Watson (see \cite{Wa}):
\begin{eqnarray}\label{wats}
\left|\int_Y\phi_i\phi_j\phi_kdv\right|^2=
\frac{\Lambda(1/2,\phi_i\otimes\phi_j\otimes\phi_k)} {
\Lambda(1,\phi_i,Ad)\Lambda(1,\phi_j,Ad)\Lambda(1,\phi_k,Ad)}\ .
\end{eqnarray}
Here the $\phi_t$ are the so-called cuspidal Hecke-Maass functions of norm one on
the Riemann surface $Y=\G\sm \uH$ arising from the full modular group $\G=SL_2(\bz)$
or from the group of units of a quaternion algebra. The functions
$\Lambda(s,\phi_i\otimes\phi_j\otimes\phi_k)$ and $\Lambda(s,\phi,Ad)$ are
appropriate {\it completed} automorphic $L$-functions.

It was first discovered by R.~Rankin and A.~Selberg that the special cases of triple
products as above give rise to automorphic $L$-functions (namely, they considered
the case where one of Maass forms is replaced by an Eisenstein series). That allowed
them to obtain analytic continuation and effective bounds for these $L$-functions
and, as an application, to obtain first non-trivial bounds for Fourier coefficients
of cusp forms towards Ramanujan conjecture. The relation (\ref{wats}) can be viewed
as a far reaching generalization of the original Rankin-Selberg formula. The
relation \eqref{wats} was motivated by the work \cite{HK} by M.~Harris and S.~Kudla
on a conjecture of H.~Jacquet.

\subsection{Results}\label{results} \ In this paper we consider
the following problem. We fix two Maass forms $\phi = \phi_\tau$ and $
\phi'= \phi_{\tau'}$ as above and consider the coefficients defined by
the triple period:
\begin{eqnarray}\label{ci} c_i=\int_Y\phi\phi'\phi_idv\ \ \
\end{eqnarray}
as the $\phi_i$ run over an orthonormal basis of Maass forms.

Thus we see from (\ref{wats}) that the estimates of the coefficients
$c_i$ are essentially equivalent to the estimates of the
corresponding $L$-functions. One would like to have a general
method of estimating the coefficients $c_i$ and similar
quantities. This problem was raised by Selberg in his celebrated
paper \cite{Se}.

The first non-trivial observation is that the coefficients $c_i$
have exponential decay in $|\lm_i|$ as $i\to\8$. Namely, as we
have shown in \cite{BR3}, it is natural to introduce normalized
coefficients
\begin{eqnarray} d_i=\g(\lm_i)|c_i|^2\ .
\end{eqnarray}
Here $\g(\lm)$ is given by an explicit rational expression in terms of the standard
Euler $\G$-function (see \cite{BR3}) and, for purely imaginary $\lm$, it has an
asymptotic   $\g(\lm)\sim\beta |\lm|^{2}\exp(\frac{\pi}{2}|\lm|)$ when $|\lm|\to\8$
with some explicit $\beta>0$. It turns out that the normalized coefficients $d_i$
have at most {\it polynomial growth} in $|\lm_i|$, and hence the coefficients $c_i$
decay exponentially. This is consistent with (\ref{wats}) and general experience
from the analytic theory of automorphic $L$-functions (see \cite{BR3}, \cite{Wa}).

In  \cite{BR3} we proved the following mean value bound
\begin{eqnarray}\label{mean-value-convex}
\sum_{|\lm_i|\leq T}d_i\leq AT^2\ ,
\end{eqnarray}
for arbitrary $T > 1$ and some effectively computable constant $A$.

According to Weyl's law the number of terms in this sum is of order $C T^2$. So this
formula says that on average the coefficients $d_i$ are bounded by some constant.

 More precisely, let us we fix an interval $I \subset \br$ around point $T$
and consider the finite set of all Maass forms $\phi_i$ with
parameter $|\lm_i|$ inside this interval. Then the average value
of coefficients $d_i$ in this set is bounded by a constant {\it
provided} the interval $I$ is long enough (i.e., of size $\approx
T$).

Note that the best individual bound which we can get from this formula is $d_i\leq
A|\lm_i|^2$. For Hecke-Maass forms this bound corresponds to the convexity bound for
the corresponding $L$-function via Watson formula (\ref{wats}).

In this paper we outline the proof of the following bound.

\begin{theorem}\label{thm}
There exist effectively computable constants $B,\ b>0$ such that, for an arbitrary
$T > 1$ we have the following bound
\begin{eqnarray}\label{sub-convex-thm}
\sum_{|\lm_i| \in I_T} d_i\leq B T^{5/3}\ ,
\end{eqnarray}
where $I_T$ is the interval of size $b T^{1/3}$ centered at $T$.
\end{theorem}

 Note that this theorem gives an individual bound $d_i \leq B |\lm_i|^{5/3}$
\   (for $|\lm_i| > 1$). Thanks to the Watson formula (\ref{wats}) and a lower bound
of H. Iwaniec \cite{I} on $L(1,\phi_{\lm_i},Ad)$ this leads to the following
 {\it   {subconvexity}} bound for the triple $L$-function (for an exact relation between
triple period and special values of $L$-functions, see \cite{Wa}).

\begin{coroll}
Let $\phi$ and $\phi'$ be fixed Hecke-Maass cusp forms. For any
$\eps>0$,
 there exists $C_\eps>0$
such that the bound
\begin{eqnarray}\label{sub} L(\shaf\ ,
\phi\otimes\phi'\otimes\phi_{\lm_i})\leq C_\eps|\lm_i|^{5/3+\eps}
\end{eqnarray}  holds for any
Hecke-Maass form $\phi_{\lm_i}$.
\end{coroll}

The convexity bound for the triple $L$-function corresponds to
\eqref{sub} with the exponent $5/3$ replaced by $2$. We refer to
\cite{IS} for a discussion of the subconvexity problem which is in
the core of modern analytic number theory. We note that the above
bound is the first subconvexity bound for an $L$-function of
degree $8$. All previous subconvexity results were obtained for
$L$-functions of degree at most $4$.

 Recently, using ergodic
theory methods, A. Venkatesh \cite{V} obtained a subconvexity
bound for the triple $L$-function in the level aspect (i.e., with
respect to a tower of congruence subgroups $\G(N)$ as $N\to\8$).

We formulate a natural
\begin{conject}
For any $\eps>0$ we have $d_i\ll|\lm_i|^{\eps}$ .
\end{conject}

For Hecke-Maass forms on congruence subgroups, this conjecture is
consistent with the Lindel\"{o}f conjecture for the triple
$L$-functions (for more details, see \cite{BR3} and \cite{Wa}).

\subsection{Remarks.}\
1) Our results can be generalized to the case of a general finite
co-volume lattice $\G \subset G$. In this case the spectral
decomposition of the Laplace-Beltrami operator on $Y = \G
\backslash \uH$ is given by a collection of eigenfunctions
$\phi_s$ (including the Eisenstein series) where the parameter $s$
runs through some set $S$ with the Plancherel measure $d\mu$; for
any function $u \in C_c^\8(Y)$ the spectral decomposition takes
the form $\int_S |<u, \phi_s>|^2 d\mu = ||u||^2_{L^2(Y)}.$

Let us fix two Maass {\it cusp forms} $\phi$ and $\phi'$ on $Y$.
For every $s \in S$ we define the parameter $\lm_s\in\bc$ and the
coefficient $d_s$ in the same way as before. In this case we can
prove the bound
$$
\int_{S_T} d_s \ d\mu \leq B T^{5/3 \ +  \eps}\ , \ {\textrm where} \
 S_T = \{s \in S \ |\  |\lm_s| \in I_T \}
 $$
 2) First results on the exact exponential decay of triple
products for a general lattice $\G$ were obtained by A. Good \cite
{Go} and P. Sarnak \cite{Sa3} using ingenious analytic
continuation of Maass form to the complexification of the Riemann
surface $Y$ (for representation-theoretic approach to this method
and generalizations, see  \cite{BR1} and \cite{KS}). Our present
method seems to be completely different and avoids analytic
continuation.
\section {The method}

We describe now the general ideas behind our proof. It is based on
ideas from representation theory (for a detailed account of the
corresponding setting, see \cite{BR3}). In what follows we sketch
the method of the proof with the complete details appearing
elsewhere.

\subsection{Automorphic representations}{\label{automorphic}
 \ Let $G$ denote the group of all motions of $\uH$. This group is naturally
 isomorphic to $\PGLR$ and as a $G$-space $\uH$ is naturally isomorphic to
  $G / K$, where $K = PO(2)$ is the standard maximal compact subgroup of $G$.

  By definition,  $\G$ is a subgroup of $G$. The space $X = \G \backslash G$ with the
  natural right action of $G$ is called an {\it automorphic space}. We will identify
  the Riemann surface $Y = \G \backslash \uH$ with $X / K$.

We start with the fact that every automorphic function $\phi$
(e.g., a Maass form) generates an automorphic representation of
the group $G$; this means that, starting from $\phi$, we produce a
smooth irreducible pre-unitary representation of the group $G$ in
a space $V$ and its isometric realization $\nu : V \to C^{\8}(X)$
in the space of smooth functions on $X$. If a Maass form $\phi$
has the eigenvalue $\mu=\frac{1-\lm^2}{4}$ then the corresponding
representation V is isomorphic to the representation of the
principal series $V_\lm$ when $\lm \in i \br$, to the
representation of complementary series $V_\lm$ when $\lm \in
[0,1)$ and to the trivial representation when $\lm = 1$.

This means that we have a very explicit {\it model} \ of the
abstract subspace $V \subset  C^\8(X)$ as the space of smooth even
homogeneous functions on $\br^2 \setminus 0$ \ of homogeneous
degree $\lm - 1$. Restricting to the unit circle $S^1 \subset
\br^2$, we get realization of $V$ as the space of smooth even
functions on the circle $S^1$ (for details, see \cite {BR3}). We
will use this model to make explicit computations.

The triple product $c_i=\int_Y\phi\phi'\phi_idv$ extends to a
$G$-equivariant trilinear functional on the corresponding
automorphic representations $l^{aut}_i:V\otimes V'\otimes
V_i\to\bc$, where $V = V_\tau, V' = V_{\tau'}$ and $V_i =
V_{\lambda_i}$ .

Next we use a general result from representation theory that such
$G$-equivariant trilinear functional is unique up to a scalar
(\cite{O}, \cite{Pr}). This implies that the automorphic
functional $l^{aut}_i$ is proportional to some explicit {\it
model} functional $l^{mod}_{\lm_i}$.  In \cite{BR3} we gave a
description of such model functional $l^{mod}_{\lm}: V \otimes V'
\otimes V_{\lm} \to \bc$ for any $\lm$ using explicit realizations
of representations $V$, $V'$ and $V_{\lm}$ of the group $G$ in
spaces of homogeneous functions; it is important that the model
functional knows nothing about automorphic picture and carries no
arithmetic information.

Thus we can write $l^{aut}_i = a_i \cdot l^{mod}_{\lambda_i}$ for some constant
$a_i$, and hence
\begin{eqnarray}\label{c-a} c_i=
l^{aut}_i(e_\tau\otimes e_{\tau'} \otimes e_{\lambda_i}) = a_i \cdot
l^{mod}_{\lm_i}(e_\tau \otimes e_{\tau'} \otimes e_{\lambda_i})\ ,
\end{eqnarray}
where $e_\tau,\ e_{\tau'},\ e_{\lambda_i}$ are K-invariant unit
vectors in representations $V, V'$ and $V_{\lambda_i}$
corresponding to the automorphic forms $\phi$, $\phi'$ and
$\phi_i$.

It turns out that the proportionality coefficient $ a_i$ in \eqref{c-a} carries an
important \lq\lq automorphic" information while the second factor carries no
arithmetic information and can be computed in terms of $\G$-functions using explicit
realizations of representations $V_\tau$, $V_{\tau'}$ and $V_{\lambda}$ (see
Appendix in \cite{BR3} where this computation is carried out). This second factor is
responsible for the exponential decay, while the first factor $a_i$ has a polynomial
behavior in parameter $\lm_i$. An explicit computation shows that
$|c_i|^2=\frac{1}{\g(\lm_i)}|a_i|^2$, and hence $d_i=|a_i|^2$ (where $\g(\lm)$ was
described in Section \ref{results}).

\subsection{Hermitian forms}{\label{hermitian} \ In order to estimate the quantities
$d_i$, we consider the space $E = V_\tau \otimes V_{\tau'}$ and use the fact that
the coefficients $d_i$ appear in the spectral decomposition of the following {\it
geometrically defined} non-negative Hermitian form $H_\Dl$ on $E$ (for a detailed
discussion, see \cite{BR3}).

Consider the space $C^\8(X\times X)$. The diagonal $\Dl:X \to X\times X$ gives rise
to the restriction morphism $r_\Dl : C^\8(X\times X)\to C^\8(X)$. We define a
non-negative Hermitian form $H_\Dl$ on $C^\8(X\times X)$ by setting $H_\Dl =
(r_\Dl)^{\ast}(P_X)$, where $P_X$ is the standard $L^2$ Hermitian form on $C^\8(X)$
i.e.,
$$H_\Dl(w)= P_X(r_\Dl(w)) = \int_{X}|r_\Dl(w)|^2d\mu_X $$
 for any $w\in C^\8(X\times X)$.
We call the restriction of the Hermitian form $H_\Dl$ to the
subspace $E\subset C^\8(X\times X)$ the {\it diagonal} Hermitian
form and denote it by the same letter.

We will  describe the spectral decomposition of the Hermitian form $H_\Dl$ in terms
of Hermitian forms corresponding to trilinear functionals. Namely, if $L$ is a
pre-unitary representation of $G$ with $G$-invariant norm $|| \ ||_L$ then every
$G$-invariant trilinear functional  $l: V \otimes V' \otimes L \to \bc$,
 defines a Hermitian form $H^l$ on $E$ by  $H^l(w)
=\sup\limits_{||u||_L = 1} |l(w \otimes u)|^2   \ $.

 Here is another description of this form (see
\cite{BR3}). Functional $l: V \otimes V' \otimes L \to \bc$ gives
rise to a $G$-intertwining morphism $T^l: E \to L^*$ which image
lies in the smooth part of $L^*$. Then the form  $H^l$ is just the
inverse image of the Hermitian form on $L^*$ corresponding to the
inner product on $L$.

Consider the orthogonal decomposition $L^2(X) =\left(\oplus_i
V_i\right)\oplus\left( \oplus_{\kappa} V_{\kappa}\right)$ where
$V_i$ correspond to Maass forms and $V_{\kappa}$ correspond to
representations of discrete series. Every subspace $L \subset
L^2(X)$ defines a trilinear functional $l :E \otimes L \to \bc$
and hence a Hermitian form $H^l$ on $E$. Hence, the decomposition
of $L^2(X)$ gives rise to the corresponding  decomposition $H_\Dl
= \sum H_i^{aut} + \sum H_{\kappa}^{aut}$ of Hermitian forms (see
\cite{BR3}).

We denote by $H_\lm$ the {\it model} Hermitian form corresponding to the {\it model}
trilinear functional $l^{mod}_\lm :V \otimes V' \otimes V_\lm \to \bc$.  From
definition we see that $H_i^{aut} = d_i H_{\lm_i}$ which leads us to

{\bf Basic identity}
\begin{eqnarray}\label{besseli} H_\Dl\ = \sum_i d_iH_{\lm_i} + \sum_{\kap}
H_{\kap}^{aut},
\end{eqnarray}
We will mostly use the fact that for every vector $w \in E$ this basic identity
gives us an inequality
\begin{eqnarray}\label{bessel}   \sum_i d_i H_{\lm_i}(w) \leq H_\Dl (w)
\end{eqnarray}
which is an equality if the vector $r_\Dl(w)$ does not have
projection on discrete series representations (for example, if the
vector $w$ is invariant with respect to the diagonal action of $K$
on $E$).

    We can use this inequality to bound coefficients $d_i$. Namely,
    for a given vector $w \in E$ we usually can compute the values
    $H_{\lm}(w)$ by explicit computations in the model of representations $V, V',
    V_{\lambda}$. It is usually much more difficult  to get reasonable estimates of the
    right hand side $H_{\Dl}(w)$. In cases when we manage to do this  we  get
    some bounds for the  coefficients $d_i$.

\subsection{Mean-value estimates}
\ In \cite{BR3}, using the geometric properties of the diagonal form and explicit
estimates of forms $H_\lm$, we established the mean-value bound
(\ref{mean-value-convex}): $\sum\limits_{|\lm_i| \leq T} d_i \leq A T^2\ .$ Roughly
speaking, the proof of this bound is based on the fact that while the value of the
form $H_{\Dl}$ on a given vector $w \in E$ is very difficult to control, we can show
that for many vectors $w$  the value $H_\Dl(w)$ can be bounded by $P_E(w)$, where
$P_E$ is the Hermitian form which defines the standard unitary structure on $E$.

More precisely, consider the natural representation $\sigma = \pi
\otimes \pi'$ of the group $G\times G$ on the space $E$. Then for
a given compact neighborhood  $U \subset G \times G$ of the
identity element, there exists a constant $C$ such that for any
vector $w \in E$, the inequality $H_{\Dl}(\sigma(g)w) \leq C
P_E(w)$ holds for at least half of the points $g \in U$. This
follows from the fact that the average over $U$ of the quantity
$H_{\Dl}(\sigma(g)w)$ is bounded by $C P_E(w)/2$.

This allows us for every $T\geq 1$, to find a vector $w \in E$ such that $H_{\Dl}(w)
\leq C T^2$ while  the inequality $H_\lm(w) \geq  c$ \ holds for all $|\lm| \leq T$.

\subsection{Bounds for sums over shorter intervals} \
   The main starting point of our approach to the
subconvexity bound is the inequality (\ref{bessel}) for Hermitian forms.  For a
given $T>1$, we construct a test vector $w_T\in E$ such that the weight function
$\lm\mapsto H_\lm(w_T)$ has a sharp peak near $|\lm|=T$ (i.e., a vector satisfying
the condition (\ref{(i)}) below).

The problem is how to estimate effectively $H_\Dl(w_T)$. The idea is that the
Hermitian form $H_\Dl$ is geometrically defined and, as a result, satisfies some
non-trivial bounds, symmetries, etc. None of the explicit {\it model} Hermitian
forms $H_\lm$ satisfies similar properties. By applying these symmetries to the
vector $w_T$ we construct a new vector $\tilde w_T$ and from the geometry of the
automorphic space $X$ we deduce the bound $H_\Dl(w_T)\leq H_\Dl(\tilde{w}_T)$.

On the other hand, the weight function $H_\lm(\tilde w_T)$ in the
spectral decomposition $H_\Dl(\tilde w_T) = \sum d_i H_{\lm_i}(
\tilde w_T)$ for $\tilde w_T$ behaves quite differently from the
weight function $H_\lm(w_T)$ for $w_T$. Namely, the function
$H_\lm(\tilde{w}_T)$ behaves regularly (i.e., satisfies condition
(\ref{(ii)}) below), while the weight function $H_\lm(w_T)$ has a
sharp peak near $|\lm| = T$.

The regularity of the function $H_\lm(\tilde w_T)$ coupled with the mean-value bound
(\ref{mean-value-convex}) allows us to prove a sharp upper bound on the value of
$H_\Dl(\tilde{w}_T)$ by purely spectral considerations (in cases we consider there
is no contribution from discrete series). We do not see how to get such sharp bound
by geometric considerations.

  Using this bound for $H_\Dl(\tilde w_T)$ and the inequality $H_\Dl(w_T)\leq H_\Dl(\tilde{w}_T)$
  we obtain a non-trivial bound
for $H_\Dl(w_T)$ and, as a result, the desired bound for the coefficients $d_i$.

\subsection{Formulas for test vectors}
 \  Let us describe the construction of vectors $w_T, \tilde w_T$.
We  assume for simplicity that $V'\simeq \bar V$ -- the complex
conjugate representation; it  is also an automorphic
representation with the realization $\bar\nu:\bar V\to C^\8(X)$.
It is easy to see that the upper bound estimate that we need in
the general case can be easily reduced to this special case.

We only consider the case of representations of the principal series, i.e. we assume
that $V = V_{\tau}$, $V' =\bar V= V_{-\tau}$ for some $ \tau \in i \mathbb{R}$; \
the case of representations of the complementary series can be treated similarly.

 Let $\{e_n\}_{n\in
2\bz}$ be a $K$-type orthonormal basis in $V$. We denote by $\{e'_n=\bar e_{-n}\}$
the complex conjugate basis in $\bar V$.

For a given  $T\geq 1$ we choose even $n$ such that $|T-2n|\leq
10$ and set
$$w_T=e_n\otimes e'_{-n}\ \ \ \text{ and}\ \ \ \ \tilde w_T=e_n\otimes
e'_{-n}+e_{n+2}\otimes e'_{-n-2}\ .$$

With such a choice of test vectors we have the following bounds.

{\bf Geometric bound:}
\begin{eqnarray}
\label{(star)} H_\Dl(w_T)\leq H_\Dl(\tilde w_T)
\end{eqnarray}

{\bf Spectral bounds:}

(i) There exist constants $b, c > 0$ such that
\begin{eqnarray}
\label{(i)} H_\lm(w_T)\geq c|\lm|^{-5/3}\text{ for } |\lm| \in I_T\
\end{eqnarray}
where $I_T$ is the interval of length $bT^{1/3}$ centered at point
$T$.

(ii) There exists a constant $c'$ such that
\begin{eqnarray}\label{(ii)}H_\lm(\tilde w_T)\leq
\begin{cases}
c'T\inv(1+|\lm|)^{-1}&\text{for all}\ |\lm|\leq 2T\ ,\\
   c'|\lm|^{-3}&\text{for all}\ |\lm|>2T\ . \\
\end{cases}
\end{eqnarray}

Using the bound \eqref{(ii)}  we can get a sharp estimate of
$H_\Dl(\tilde w)$. Namely, from \eqref{besseli} we conclude that
  $H_\Dl(\tilde w) = \sum d_i H_{\lm_i}(\tilde w)$
  (since vectors $\tilde w_T$ are $\Dl K$-invariant, we do not
have  contribution  from representations of discrete series).

  The spectral bound (\ref{(ii)}) for $H_\lm(\tilde w)$  together with
   the mean-value bound
  (\ref{mean-value-convex}) for coefficients $d_i$ imply that
$$
H_\Dl(\tilde w_T)\leq D
$$ for some explicit constant $D$.

 Using the geometric
inequality (\ref{(star)}) we see that $H_\Dl(w_T)\leq D$. Using the spectral bound
(\ref{(i)}) we obtain
$$
\sum_{|\lm_i| \in I_T}d_icT^{-5/3}\leq \sum_{i}d_iH_{\lm_i}(w_T)\leq H_\Dl(w_T)\leq
D\ .
 $$
From this we deduce the bound (\ref{sub-convex-thm}) in Theorem
\ref{thm}.

\subsection{Proof of the geometric bound {\ref{(star)}}.} \ The
inequality (\ref{(star)}) easily follows from the pointwise bound
on $X$ due to the fact that, in the automorphic realization, the
vector $e_n\otimes e'_{-n}$ is represented by a function which
restriction $u_n = r_{\Delta}(e_n \otimes e'_{-n})$ to the
diagonal is non-negative
$$
u_n(x) = \nu(e_n)(x)\cdot\bar\nu(e'_{-n})(x)=|\nu(e_n)(x)|^2\geq 0.
$$
From this we see that
$
H_\Dl(w_T)=\int_{X}|u_n(x)|^2d\mu_X \leq
\int_{X}|u_n(x)+u_{n+2}(x)|^2d\mu_X= H_\Dl(\tilde{w}_T)\ .
$

\subsection{Sketch of proof of the spectral bounds (\ref{(i)}) and
(\ref{(ii)}).}\ We will use the explicit form of the kernel defining Hermitian forms
$H_\lm$ in the model realizations of representations $V$, $V'$ and $V_\lm$. Namely,
we use the standard realization of these representations in the space
$C^\8_{even}(S^1)$ of even functions on $S^1$ (see \cite{BR3} and
\ref{automorphic}). Under this identification the basis $\{e_n\}$ becomes the
standard basis of exponents $\{e_n=e^{in\theta}\}$, where $0\leq\theta<2\pi$ is the
standard parameter on $S^1$.

 As was
shown in \cite{BR3}, Section 5, in such realization the invariant functional
$l^{mod}_\lm$ on the space $V\otimes  V' \otimes V_\lm\simeq C^\8((S^1)^3)$ is given
by the following kernel on $(S^1)^3$
$$L_\lm(\theta,\theta',\theta'')=
|\sin(\theta-\theta')|^{\frac{-1+\lm}{2}}|\sin
(\theta-\theta'')|^{\frac{-1+2\tau-\lm}{2}}|\sin
(\theta'-\theta'')|^{\frac{-1-2\tau-\lm}{2}}\ ,$$ where $V = V_\tau$, $ V' =
V_{-\tau}$ with $\tau\in i\br$. From this it follows that the Hermitian forms
$H_\lm$ on $E\simeq C^\8(S^1\times S^1) $ are given by oscillatory integrals (over
$(S^1)^4$) and the verification of conditions (\ref{(i)}) and (\ref{(ii)}) is
reduced to the stationary phase method.

In fact we will use the values of $H_\lm(w)$ only for $\Dl
K$-invariant vectors $w \in E$. This considerably simplifies our
computations since we can reduce them to two repeated integrations
in one variable and use the stationary phase method in {\it one
variable}.

Namely, let us fix a $\Dl K$-invariant vector $w \in E$. Then the
vector $T_\lm(w)\in V_{-\lm}$ is proportional to the standard
$K$-invariant vector $e_0\in V_{-\lm}\simeq C^\8_{even}(S^1)$
(here the operator $T_\lm : E \to V_{-\lm}$ corresponds to the
model trilinear functional $l^{mod}_\lm$ as described in
\ref{hermitian}). This implies that $H_\lm(w) = |T_\lm(w)(0)|^2$.
The value $T_\lm(w)(0)$ is given by the following oscillating
integral
$$T_\lm(w)(0) = <w,K_\lm>=\int w(\theta,\theta')K_\lm(\theta,\theta')d\theta d\theta'
,$$ where $ K_\lm(\theta,\theta')=L_\lm(\theta,\theta',0). $ Since
the vector $w$ is $\Dl K$-invariant it can be described by a
function  in one variable; namely, $w(\theta,\theta')= u(c)$ for
$u \in C^\8(S^1)$ and $c = (\theta - \theta')/2$. We have $<w,
K_\lm> = \int u(c)k_\lm(c)dc$, where the function $k_\lm$ is
obtained from $K_\lm$ by averaging over $\Dl K$. Thus for a $\Dl
K$-invariant vector $w$, the estimates of $H_\lm(w)$ are
equivalent to estimates of the one-dimensional integral $<u,
k_\lm> = \int u(c) k_\lm(c)dc$.

The function $k_\lm(c)$, which is obtained from
$K_\lm(\theta,\theta')$ via one-dimensional integration, is not an
elementary function. However, using stationary phase method, we
obtain the representation
$k_\lm(c)=|\lm|^{-\haf}m_\lm(c)+r_\lm(c)$, where the main term
$m_\lm$ (given by contributions from non-degenerate stationary
points of the phase in the corresponding integral) is  an
elementary function
$$ m_\lm(c)=\al(\lm)|\sin(c)|^{-\haf-\frac{\lm}{2}}|\cos(c)|^{-\haf+\frac{\lm}{2}}\
\ ,$$ with $\al(\lm)=(\pi)\inv
e^{-i\frac{\pi}{4}}2^{-\haf+\frac{\lm}{2}}$. The stationary phase
method also gives a  bound for  the remainder term
$||r_\lm||_{L^1(S^1)}\leq a (1+|\lm|)^{-3/2}$ for some constant
$a$.

The vectors $w$ which we consider correspond to bounded functions $u(c)$. For such
vectors, the estimate of $<u, k_\lm>$ is reduced to the estimate of $<u, m_\lm>\ =
\int u(c) m_\lm(c)dc$.

We deduce spectral bounds (\ref{(i)}) and (\ref{(ii)}) by applying
stationary phase method to integrals $<u_T, m_\lm>$ and $<\tilde
u_T,m_\lm>$, where $u_T,\ \tilde u_T \in C^\8(S^1)$ are functions
corresponding to vectors $w_T,\ \tilde w_T \in E$.

The key fact responsible for the crucial bound in (\ref{(i)}) is that for $T=|\lm|$,
the phase of the oscillating integral $<u_T,m_\lm>$ corresponding to the value
$H_\lm(w_T)$ has a {\it degenerate} critical point at $c = \pi /4$ with the {\it
non-vanishing} amplitude at that point. For other values of $\lm$ this phase has
nondegenerate critical points.

 Since this critical point is degenerate the integral
$<u_T,m_\lm>$ has a sharp peak at $|\lm| = T$. The standard technique developed to
analyze the asymptotic behavior of the Airy functions then gives the bound
(\ref{(i)}) for $|\lm| \asymp T $.

On the other hand for the oscillating integral $<\tilde u_T,m_\lm>$ corresponding to
the value $H_\lm(\tilde w_T)$ the phase is the same as for the integral $<u_T,
m_\lm>$,\  but the amplitude has an additional factor $a(c) = 1 + e^{4ic}$ which was
chosen in such a way that it vanishes at the degenerate critical point which
develops at $|\lm| = T$. As a result this point does not give an additional
contribution to this integral.

   This is a classical situation for which the uniform bounds for the oscillating
integrals are well-known (e.g., bounds on the Airy function and its derivative).
From this we deduce the bound in (\ref{(ii)}). In fact, we find that for $|\lm|>T$
there are no critical points at all. This implies that for $|\lm|> 2T$ we have a
stronger bound $H_\lm(\tilde w_T)\ll |\lm|^{-N}$  for any $N>1$ (compare to
\eqref{(ii)}).

For $|\lm|\ll T$, we also consider singularities of the amplitude in the
corresponding integrals in order to show that the low-lying spectrum contribution is
bounded. This includes the contribution from representations of the complementary
series and the trivial representation (in fact, we have to deal with the
singularities of the amplitude for all values of $\lm$).

The above arguments also  prove the following result on the $L^4$-norm of $K$-types
in irreducible automorphic representations of $\PGLR$. This result is of
independent interest.

\begin{theorem} For a fixed class
one automorphic representation $\nu:V\to C^\8(X)$, there exists
$D>0$ such that $||\nu(e_n)||_{L^4(X)}\leq D$ for all $n$.
\end{theorem}

One would expect that a similar fact holds for representations
of the discrete series as well.

\subsection*{Acknowledgements}\
It is a pleasure to thank Peter Sarnak for stimulating discussions
and support of this work. We would like to thank Herv\'{e} Jacquet
for a valuable comment.

Research was partially supported by BSF grant, by Minerva
Foundation and by the Excellency Center ``Group Theoretic Methods
in the Study of Algebraic Varieties'' of the Israel Science
Foundation, the Emmy Noether Institute for Mathematics (the Center
of Minerva Foundation of Germany), and by RTN ``Liegrits''. The
results of this paper were obtained during our visits to
Max-Planck Institute in Bonn and in Leipzig and to Courant
Institute.

\bibliographystyle{amsplain}

\end{document}